\documentclass[]{cmip-l}

 \usepackage[francais]{babel} 

\def\abstract{\if@twocolumn
\section*{R\'esum\'e}
\else \small
\begin{center}
{\bf R\'esum\'e\vspace{-.5em}\vspace{0pt}}
\end{center}
\quotation \fi}

\def\today{
 \number\day\space
 \ifcase\month\or
 janvier\or {f\'evrier}\or mars\or avril\or mai\or juin\or
 juillet\or {ao\^ut}\or septembre\or octobre\or novembre
 \or {d\'ecembre}\fi
 \space\number\year}

\def\tableofcontents{
 \section*{Table des mati\`eres\markboth{{Table des mati\`eres}}{{Table des mati\`eres}}}
 \@starttoc{toc}}

\def\listoffigures{\section*{Liste des figures\markboth{Liste des figures}{Liste des figures}}
  \@starttoc{lof}}

\def\thebibliography#1{\section*{Bibliographie\markboth{Bibliographie}{Bibliographie}}\list
 {[\arabic{enumi}]}{\settowidth\labelwidth{[#1]}\leftmargin\labelwidth
 \advance\leftmargin\labelsep
 \usecounter{enumi}}
 \def\newblock{\hskip .11em plus .33em minus -.07em}
 \sloppy
 \sfcode`\.=1000\relax}

\def\fnum@table{Tableau \thetable}

\def\@part[#1]#2{\ifnum \c@secnumdepth >\m@ne \refstepcounter{part}
\addcontentsline{toc}{part}{\thepart \hspace{1em}#1}\else
\addcontentsline{toc}{part}{#1}\fi { \parindent 0pt \raggedright
 \ifnum \c@secnumdepth >\m@ne \Large \bf Partie \thepart \par \nobreak \fi \huge
\bf #2\markboth{}{}\par } \nobreak \vskip 3ex \@afterheading }

\newtheorem{theo}{Th\'eor\`eme}[section]

\newtheorem{prop}[theo]{Proposition}
\newtheorem{lem}[theo]{Lemme}
\newtheorem{rema}[theo]{Remarque}
\newtheorem{remas}[theo]{Remarques}
\newtheorem{cor}[theo]{Corollaire}

\def \kbar {{\bar k}}

\def \Xbar {{\overline X}}

\def \Romannumeral #1 {\expandafter\uppercase\expandafter {\romannumeral #1} }

\def \spec {{\rm{Spec\,}}}
\def \dim {{\rm{dim\,}}}

\def \Z {{\bf Z}}
\def \Q {{\bf Q}}
\def \F {{\bf F}}

\def \PP {{\bf P}}

\def \im {{\rm {Im\,}}}
\def \G {{\bf G}_m}
\def \gal {{\rm Gal}\,}
\def \br {{\rm Br}\,}
\def \pic {{\rm Pic}\,}

\def\smallsquare{\vbox{\hrule\hbox{\vrule height 1 ex\kern 1 ex\vrule}\hrule}}
\def\enddem{\hfill \smallsquare\vskip 3mm}

\usepackage{amscd}
\usepackage{amssymb}
\usepackage{amsmath}
\usepackage{pb-diagram}

\everymath{\displaystyle}

\def \et {\hbox{\scriptsize \'et}}

\def \id {{\rm id}}

\DeclareFontFamily{U}{wncy}{}
\DeclareFontShape{U}{wncy}{m}{n}{%
   <5>wncyr5%
   <6>wncyr6%
   <7>wncyr7%
   <8>wncyr8%
   <9>wncyr9%
   <10>wncyr10%
   <11>wncyr10%
   <12>wncyr6%
   <14>wncyr7%
   <17>wncyr8%
   <20>wncyr10%
   <25>wncyr10}{}
\DeclareMathAlphabet{\cyrille}{U}{wncy}{m}{n}

\def \C{{\bf C}}

\title{Autour de la conjecture de Tate \mbox{\`a coefficients $\Z_\ell$} pour les vari\'et\'es sur les corps finis}
\author{Jean-Louis Colliot-Th\'el\`ene et Tam\'as Szamuely}
\address{C.N.R.S., U.M.R. 8628, Universit\'e de Paris-Sud, Math\'ematique, B\^atiment 425, 91405 Orsay, France}
\email{jlct@math.u-psud.fr}\date{\today}
\address{Alfr\'ed R\'enyi Institute of Mathematics, Hungarian Academy of Sciences, Re\'altanoda utca 13--15, H-1053 Budapest, Hungary}
\begin{document}
\email{szamuely@renyi.hu}\maketitle \markboth{Jean-Louis Colliot-Th\'el\`ene et Tam\'as Szamuely}{Autour de la
conjecture de Tate \`a coefficients $\Z_\ell$}
\section{Introduction}

Soient $k$ un corps fini, $\kbar$ une cl\^oture alg\'ebrique de $k$,
$G$ le groupe de Galois $\gal(\kbar|k)$ et $\ell$ un nombre premier inversible dans $k$.
Consid\'erons  une vari\'et\'e projective, lisse, g\'eom\'etriquement int\`egre $X$, de
dimension $d$. D'apr\`es la conjecture de Tate, l'application cycle \`a valeurs dans la
cohomologie \'etale $\ell$-adique induit une {\em surjection}
\begin{equation}\label{tate1}
CH^i(X)\otimes_\Z\Q_\ell\twoheadrightarrow H^{2i}(\Xbar, \Q_\ell(i))^G.
\end{equation}
Une forme \'equivalente de la conjecture est la surjectivit\'e du morphisme
\begin{equation}\label{tate2}
CH^i(\Xbar)\otimes_\Z\Q_\ell\to\bigcup_U H^{2i}(\Xbar, \Q_\ell(i))^U
\end{equation}
o\`u $\Xbar:=X\times_k\kbar$ et $U$ parcourt le syst\`eme des sous-groupes ouverts de $G$.
La forme plus forte ci-dessus en r\'esulte par un argument de restriction-corestriction.

On peut \'egalement consid\'erer des formes enti\`eres de ces \'enonc\'es, et se demander
si les morphismes

\begin{equation}\label{tateZ1}
CH^i(X)\otimes_\Z\Z_\ell\to H^{2i}(\Xbar, \Z_\ell(i))^G
\end{equation}
ou
\begin{equation}\label{tateZ2}
CH^i(\Xbar)\otimes_\Z\Z_\ell\to\bigcup_U H^{2i}(\Xbar, \Z_\ell(i))^U
\end{equation}
induits par l'application cycle sont surjectifs. Ici le deuxi\`eme \'enonc\'e de
surjectivit\'e est a priori plus faible.

Comme nous allons le rappeler dans la section \ref{conjfausse}, on ne s'attend pas \`a ce
que les formes enti\`eres de la conjecture ci-dessus soient vraies.
N\'eanmoins, il est raisonnable d'esp\'erer
la surjectivit\'e de (\ref{tateZ1}) et (\ref{tateZ2}) pour $i=d-1$, i.e. pour les
1-cycles.

Dans ce cas, la surjectivit\'e de (\ref{tateZ2}) a \'et\'e conditionnellement
d\'emontr\'ee par Chad Schoen   :

\begin{theo}\label{theoschoen} {\rm (Schoen \cite{schoen})}
Soient $k$, $G$ et $X$ comme ci-dessus. Supposons la conjecture de Tate connue
pour les diviseurs sur une surface projective et  lisse  sur un corps fini. Alors le morphisme
$$
CH_1(\Xbar)\otimes\Z_\ell\to \bigcup_UH^{2d-2}(\Xbar, \Z_\ell(d-1))^U
$$
est surjectif, o\`u $U$ parcourt le syst\`eme des sous-groupes ouverts de  $G$.
\end{theo}

Notons que la conjecture de Tate pour les diviseurs sur une surface au-dessus d'un corps
fini peut \^etre per\c cue comme un analogue de la finitude hypoth\'etique du groupe de
Tate-Shafarevich de la jacobienne d'une courbe sur un corps de nombres.

Nous expliquons la d\'emonstration de Schoen (avec quelques modifications) dans les
sections \ref{seclef}, \ref{secalglin} et \ref{schoenfin}.

Au paragraphe \ref{brauermanin}, on voit que le th\'eor\`eme de Schoen a des cons\'e\-quences
sur l'existence de z\'ero-cycles sur certaines vari\'et\'es d\'efinies sur
un corps de fonctions d'une variable sur
la cl\^oture alg\'ebrique
d'un corps fini. Voici un cas particulier concret  :

\begin{cor}\label{corschoen} Soient  $\kbar$ et $X$ comme ci-dessus. Supposons qu'il exis\-te un
\mbox{$\kbar$-morphisme} propre surjectif  $f : \Xbar \to \overline C$, avec $\overline C$
une $\kbar$-courbe propre lisse. Supposons en outre que la fibre g\'en\'erique de $f$ est
une intersection compl\`ete lisse de dimension $\geq 3$ et de degr\'e premier \`a ${\rm
car}(k)$ dans un espace projectif, et que chacune des fibres de $f$ poss\`ede une
composante de multiplicit\'e 1. Si la conjecture de Tate pour les diviseurs sur les
surfaces projectives lisses sur un corps fini est vraie, alors le pgcd des degr\'es des
multisections de $\Xbar \to\overline C$ est \'egal \`a 1. \end{cor}

\section{G\'en\'eralit\'es sur la conjecture de Tate \`a coefficients entiers }
\label{conjfausse}

On entend souvent dire  : la conjecture de Hodge \`a coefficients entiers est fausse,
il n'est pas raisonnable d'\'enoncer  la conjecture de Tate avec des coefficients entiers.
Quelle est la situation ?

Chacune de ces conjectures porte sur l'image  d'une application cycle \'emanant du groupe
de Chow $CH^r(X)$ des cycles de codimension $r$  sur une vari\'et\'e projective et lisse
$X$ de dimension $d$,  \`a valeurs dans un groupe de cohomologie. Il s'agit de
$H^{2r}(X,\Z)$ pour Hodge et de $ H^{2r}_{\et}(X \times_{k}{\overline
  k},\Z_{\ell}(r))$ pour Tate (dans cette section on va distinguer
les groupes de cohomologie \'etale des groupes de cohomologie singuli\`ere par des indices
pour ne pas induire une confusion). On trouvera dans le survol \cite{voisin} de Voisin un
\'etat des lieux pour la conjecture de Hodge.

Si l'on croit \`a ces conjectures \`a coefficients rationnels, la variante
enti\`ere peut \^etre mise en d\'efaut de deux fa\c cons:
\smallskip

\noindent $(a)$ on trouve une
classe de cohomologie de torsion qui n'est pas la classe d'un cycle; \smallskip

\noindent $(b)$
on trouve une classe de cohomologie d'ordre infini qui n'est pas dans l'image de
l'application cycle, mais qui donne un \'el\'ement de torsion dans son conoyau.  \smallskip

Pour la conjecture de Hodge enti\`ere,  il y a des contre-exemples de type $(a)$ dus \`a
Atiyah et Hirzebruch \cite{ah}, reconsid\'er\'es plus r\'ecemment par Totaro
(\cite{totaroJAMS}, \cite{totaro}), pour les groupes $H^{2r}(X,\Z)$ avec $r\geq 2$.
L'exemple de dimension minimale chez eux est une vari\'et\'e de dimension $7$, avec une
classe de torsion
 dans $H^4(X,\Z)$.
Dans la litt\'erature (par exemple dans Milne \cite{milne}, Aside 1.4) il est affirm\'e
que l'on peut adapter ces exemples  pour donner des contre-exemples \`a la conjecture de
Tate enti\`ere sous la forme (\ref{tateZ2}), mais \`a notre connaissance aucune
d\'emons\-tration n'a \'et\'e \'ecrite. Voici donc une esquisse de d\'emons\-tration qui met
en relief les modifications \`a faire par rapport au cas analytique discut\'e dans
\cite{ah}. Il s'agit de prouver le th\'eor\`eme suivant :

\begin{theo}\label{AH} ${}$
\begin{enumerate}
\item Soit $V$ une vari\'et\'e projective et lisse sur un corps alg\'ebriquement clos. Pour
tout $\ell\geq i$ inversible sur $V$ les op\'erations de Steenrod de degr\'e impair
s'annulent sur la classe de tout cycle alg\'ebrique dans le groupe $H^{2i}_{\et}(V,
\Z/\ell\Z(i))$.
\item Au-dessus de tout corps alg\'ebriquement clos,
pour tout premier $\ell$ diff\'erent de la caract\'eristique, il existe une intersection
compl\`ete lisse $Y\subset \PP^N$, un groupe fini $G$ agissant  librement sur $Y$, une
classe $c$ de $\ell$-torsion dans $H^{4}_{\et}(Y/G, \Z_\ell(2))$ et une op\'eration  de
Steenrod de degr\'e impair qui n'annule pas
 l'image de  $c$
dans $H^{4}_{\et}(Y/G, \Z/\ell\Z(2))$.
\end{enumerate}
\end{theo}

Ici pour $\ell>2$ premier
\og
op\'eration de Steenrod de degr\'e impair
\fg{}
veut dire un
compos\'e d'op\'erations de Steenrod ${\mathcal P}^j$ et d'une op\'eration de Bockstein.
Les op\'erations ${\mathcal P}^j:\, H^{i}_{\et}(V, \Z/\ell\Z)\to H^{i+2j(\ell-1)}_{\et}(V,
\Z/\ell\Z)$ en cohomologie \'etale ont \'et\'e d\'efinies par Mme Raynaud dans \cite{ray}.
Pour $\ell=2$ on utilise des op\'erations ${Sq}^j:\, H^{i}_{\et}(V, \Z/2\Z)\to
H^{i+j}_{\et}(V, \Z/2\Z)$, \'egalement d\'efinies dans \cite{ray}.

Si le corps de base est une cl\^oture alg\'ebrique $\overline F$ d'un sous-corps $F$,
toute classe de torsion dans $H^{2i}_{\et}(V, \Z_\ell(i))$ est invariante par un
sous-groupe ouvert de $\gal(\overline F|F)$, donc pour $F$ fini le th\'eor\`eme nous
fournit un contre-exemple du type $(a)$ \`a la surjectivit\'e des applications
(\ref{tateZ1}) et (\ref{tateZ2}).

\smallskip

Esquissons une d\'emonstration du th\'eor\`eme qui nous a \'et\'e g\'en\'ereuse\-ment
communiqu\'ee par Burt Totaro. Pour d\'emontrer (1), la premi\`ere observation est que par
le th\'eor\`eme de Riemann-Roch sans d\'enomina\-teurs de Jouanolou (\cite{fulton},
Example 15.3.6) pour $\ell$ premier \`a $(i-1)!$ toute classe de cycle dans
$H^{2i}_{\et}(V, \Z/\ell\Z(i))$ est combinaison lin\'eaire de classes de Chern $c_i(E)$ de
fibr\'es vectoriels $E$ sur $X$. Il suffit donc de d\'emontrer l'\'enonc\'e d'annulation
pour les $c_i(E)$. Un calcul d'op\'erations de Steenrod montre que l'annulation vaut pour
$E$ si et seulement si elle vaut pour $E\otimes L$ avec $L$ tr\`es ample de rang un.
Ainsi, on peut supposer que $E$ est engendr\'e par ses sections globales, et a fortiori
qu'il est la tirette du fibr\'e tautologique d'une grassmannienne. L'\'enonc\'e r\'esulte
alors de la fonctorialit\'e contravariante des ${\mathcal P}^j$
et de la trivialit\'e de
la cohomologie d'une grassmannienne en degr\'es impairs (\cite{sga5}, expos\'e VII,
proposition 5.2).

Un point clef de l'argument d'Atiyah--Hirzebruch \cite{ah} \'etait l'identi\-fication de
la cohomologie en bas degr\'es d'une vari\'et\'e de Godeaux--Serre $Y/G$ comme dans (2)
\`a celle du produit d'espaces classifiants $BG\times B\G$.  On peut alg\'ebriser leur
m\'ethode en utilisant l'approxima\-tion alg\'ebrique de $BG$ introduite par Totaro. En
effet, d'apr\`es (\cite{totaro}, Remark 1.4) pour tout $s\geq 0$ il existe une
repr\'esentation $k$-lin\'eaire $W$ de $G$ telle que l'action de $G$ soit libre en dehors
d'un ferm\'e $S$ de codimension $s$ dans $W$. La cohomologie de $BG:=(W\setminus S)/G$ est
\'egale \`a celle de $G$ jusqu'en degr\'e $s\,{}$; en particulier, elle ne d\'epend
ni du choix de $W$ ni du choix de  $S$.
Le quotient $\PP(W)//G:=(\PP(W)\times (W\setminus S))/G$ est un fibr\'e
projectif sur $BG$, donc son anneau de cohomologie est un anneau de polyn\^omes sur celui
de $BG$. En particulier, la cohomologie de $BG$ est facteur direct dans celle de
$\PP(W)//G$.

Si maintenant $Y\subset \PP(W)$ est une intersection compl\`ete lisse sur laquelle
\mbox{l'action} de $G$ est libre, la cohomologie de $Y/G$ est isomorphe \`a celle de
$\PP(W)//G$ en bas degr\'es. En effet, la cohomologie de $Y$ est isomorphe \`a celle de
$\PP(W)$ jusqu'en degr\'e
$\dim(Y)-1$ par le th\'eor\`eme de Lefschetz faible. On en d\'eduit un isomorphisme entre
les cohomologies de $(Y\times (W\setminus S))/G$ et de $\PP(W)//G$  jusqu'en degr\'e
$\dim(Y)$ en appliquant la suite spectrale de Hochschild--Serre aux $G$-rev\^etements
${Y\times (W\setminus S)}\to (Y\times (W\setminus S))/G$ et $(\PP(W)\times (W\setminus
S))\to \PP(W)//G$. Or la cohomologie de $(Y\times (W\setminus S))/G$ s'identifie \`a celle
de $(Y\times W)/G$ jusqu'en degr\'e $s$, et finalement \`a celle de $Y/G$ dans
le m\^eme intervalle, puisque $W$ est un espace affine.

En somme, en bas degr\'es la cohomologie de $BG$ (donc celle de $G$) s'identifie \`a un
facteur direct de celle de $Y/G$ ci-dessus. La fin de la d\'emonstration de (2) est alors
similaire \`a celle de (\cite{ah}, Proposition 6.7). Prenons $G=(\Z/\ell\Z)^3$. Comme $G$
est d'exposant $\ell$, la suite exacte longue associ\'ee \`a ${0\to\Z_\ell\to\Z_\ell\to
\Z/\ell\Z\to 0}$
montre que $H^i(G, \Z_\ell)$ s'identifie au noyau du morphisme de Bockstein $\beta:\,
H^i(G, \Z/\ell\Z)\to H^{i+1}(G,\Z/\ell\Z)$. Le cup-produit des \'el\'ements d'une base du
$(\Z/\ell\Z)$-espace vectoriel $H^1(G,\Z/\ell\Z)\cong (\Z/\ell\Z)^3$ donne une classe dans
$H^3(G, \Z/\ell\Z)$. Essentiellement le m\^eme calcul que dans \cite{ah} montre que pour
$\ell>2$ l'image de cette classe dans $H^4(G, \Z/\ell\Z)$ par le Bockstein $\beta$ n'est
pas annul\'ee par l'op\'eration $\beta{\mathcal P}^1$, dont le degr\'e est $2\ell-1$. Pour
$\ell=2$ la m\^eme conclusion vaut pour $Sq^3$.

\medskip

Terminons cette section par une br\`eve discussion des contre-exemples de type~$(b)$. Un
c\'el\`ebre contre-exemple de ce type \`a la conjecture de Hodge a \'et\'e fabriqu\'e par
J. Koll\'ar \cite{kollar}; voir aussi \cite{soulevoisin}. Il s'agit d'une hypersurface
\og tr\`es g\'en\'erale \fg{} dans $\PP^4_{\C}$ de degr\'e $m$ un multiple de $\ell^3$ avec
$\ell$ entier premier \`a $6$, et de l'application cycle $CH^2(X) \to H^4(X,\Z)$. Comme il
s'agit d'une hypersurface de degr\'e $m$, ici on a $H^4(X,\Z)\cong\Z$, et l'image de
l'application cycle contient $m\Z$. Mais Koll\'ar montre par un argument de d\'eformation
astucieux que toute courbe sur $X$ a un degr\'e divisible par $\ell$. En d'autres mots,
l'image de $CH^2(X) \to H^4(X,\Z)$ est contenue dans $\ell H^4(X,\Z)$ et ne peut \^etre
surjective. Comme le note C. Voisin (\cite{soulevoisin}, \cite{voisin}), on peut \`a
partir de cet exemple fabriquer des contre-exemples \`a la conjecture de Hodge enti\`ere
en d'autres degr\'es aussi, par \'eclatement ou par produit direct avec une autre
vari\'et\'e.

L'\'enonc\'e de Koll\'ar en induit un au niveau de la cohomologie $\ell$-adique \'etale.
En effet, si on travaille sur un corps
alg\'ebriquement clos
non d\'enombrable et on choisit $\ell$ premier \`a
la caract\'eristique et \`a 6, sa m\'ethode fournit toujours une hypersurface
$X\subset\PP^4$ sur laquelle toute courbe a un degr\'e divisible par $\ell$. (Le corps non
d\'enombrable sert ici pour pouvoir choisir le point correspondant \`a $X$ d'un sch\'ema
de Hilbert convenable en dehors de la r\'eunion d'une famille d\'enombrable de ferm\'es
propres.) Ensuite, comme pour toute vari\'et\'e digne de ce nom, on trouve un corps $K$ de
type fini sur le corps premier sur lequel $X$ est d\'efinie. Notant $\overline K$ une
cl\^oture alg\'ebrique de $K$, l'image de l'application cycle
$$CH^2(X \times_{K}{\overline K}) \to H^4_{\et}(X \times_{K}{\overline K},\Z_{\ell}(2))\cong\Z_\ell$$
est alors contenue dans $\ell\Z_\ell$; noter qu'ici l'action de Galois sur la cohomologie
induit l'action triviale sur $\Z_\ell$.

\begin{remas}\label{remconjfausse}\rm ${}$\smallskip

\noindent 1. La m\'ethode ci-dessus ne permet pas de trouver un tel exemple  avec $K$ un
corps de nombres.\smallskip

\noindent 2. Par le th\'eor\`eme \ref{theoschoen}, si on croit \`a la conjecture de Tate
rationnelle pour les diviseurs sur les surfaces, en caract\'eristique positive le corps
$K$ ci-dessus ne peut \^etre un corps fini.\smallskip

\noindent 3. Nous ne savons pas s'il existe des contre-exemples du type $(b)$ \`a la
surjectivit\'e de (\ref{tateZ1}) sur un corps fini. En d'autres mots, nous ne savons pas
si pour $X$ projective, lisse, g\'eom\'etri\-quement connexe sur un corps fini $k$,
l'application
$$
CH^i(X)\otimes_\Z\Z_\ell\to H^{2i}_{\et}(\Xbar, \Z_\ell(i))^G/{\rm torsion}
$$
induite par l'application cycle est toujours surjective.

Cette question est \'equivalente \`a la question suivante, fort int\'eressante du point de
vue de \cite{CT} : pour tout $i\geq 0$, l'application
$$
CH^i(X)\otimes_\Z\Z_\ell\to H^{2i}_{\et}(X, \Z_\ell(i))/{\rm torsion}
$$
induite par l'application cycle
\begin{equation}\label{appcyc}
CH^i(X)\otimes_\Z\Z_\ell\to H^{2i}_{\et}(X, \Z_\ell(i)) \end{equation}
 est-elle surjective
? Le lien entre les deux questions est fourni par les suites exactes
$$ 0 \to H^1(k,H^{2i-1}_{\et}(X_{\kbar},\Z_\ell(i)))) \to H^{2i}_{\et}(X, \Z_\ell(i)) \to H^{2i}_{\et}(X_\kbar, \Z_\ell(i))^G
\to 0,$$
o\`u  les groupes $H^1(k,H^{2i-1}_{\et}(X_{\kbar},\Z_\ell(i))))$ sont des groupes finis
(ceci est une cons\'equence du th\'eor\`eme de Deligne \'etablissant les conjectures de Weil).

Notons ici pour usage ult\'erieur que par un argument bien connu utilisant  la suite de Kummer
et le groupe de Brauer,
 pour $i=1$ la surjectivit\'e de (\ref{appcyc}) est \'equivalente \`a la conjecture
de Tate \`a coefficients $\Q_\ell$, et m\^eme \`a la bijectivit\'e du morphisme
(\ref{tate1}) (voir \cite{tate}, Proposition 4.3). En vertu de la suite exacte ci-dessus,
dans le cas des diviseurs la conjecture de Tate \`a coefficients $\Q_\ell$ implique donc
la version enti\`ere sous toutes ses formes possibles.
\end{remas}

\section{Le th\'eor\`eme de Schoen, I   : un argument de type Lefschetz}\label{seclef}

Nous commen\c cons maintenant l'exposition de la d\'emonstration du th\'eor\`eme
\ref{theoschoen}, suivant \cite{schoen}.

Au cours de la preuve nous ferons \`a plusieurs reprises des extensions du corps de base
de degr\'e premier \`a $\ell$. Un argument de
restriction-corestriction fournit alors le r\'esultat au-dessus du corps de base initial.

\begin{lem} Il suffit d'\'etablir le th\'eor\`eme \ref{theoschoen} pour $d=3$.
\end{lem}

\begin{demo} D'apr\`es le th\'eor\`eme de  Bertini sur un corps fini \cite{gabber,poonen}, on peut
trouver un plongement projectif de $X$ et une hypersurface $H$ tels que le $k$-sch\'ema
$Y=X \cap H$ soit de codimension 1  dans $X$, lisse et g\'eom\'etriquement connexe.
 Comme $X\setminus Y$ est affine, pour $d>3$, les th\'eor\`emes sur la dimension
 cohomologique des sch\'emas affines
(\cite{milneEC}, \S VI.7) donnent
$$H^{2d-3}(\Xbar\setminus Y_\kbar, \Z_\ell(d-1))=0, \hskip2mm H^{2d-2}(\Xbar\setminus Y_\kbar,
\Z_\ell(d-1))=0.$$
 Donc  le morphisme compos\'e
$$
H^{2d-4}(Y_\kbar, \Z_\ell(d-2))\stackrel\sim\to H^{2d-2}_{Y_\kbar}(\Xbar, \Z_\ell(d-1))\to
H^{2d-2}(\Xbar, \Z_\ell(d-1))
$$
de l'isomorphisme de puret\'e
 et de la  fl\`eche provenant
 de la suite de localisation  est un isomorphisme  (th\'eor\`eme de Lefschetz faible, {\it ibidem}).
Comme l'application cycle est compatible aux morphismes de Gysin (\cite{milneEC},
Proposition VI.9.3), par r\'ecurrence
  sur $d$ on se ram\`ene donc au cas $d=3$.
  \end{demo}
\enddem

Jusqu'\`a la fin du paragraphe 4, on suppose donc  $d={\rm dim}(X)=3$.

\smallskip

Remarquons ensuite qu'il suffit d'\'etablir le th\'eor\`eme apr\`es avoir \'eclat\'e un
point de $X$, puisque la cohomologie de $X$ s'identifie \`a un facteur direct de celle de
l'\'eclat\'e (\cite{sga7}, expos\'e XVIII, 2.2.2) et l'application cycle est compatible
aux morphismes propres de vari\'et\'es propres et lisses (\cite{laumon}, th\'eor\`eme 6.1
et remarque 6.4).
 Ainsi, apr\`es avoir fait un \'eclatement convenable, on peut
tranquillement supposer que le deuxi\`eme nombre de Betti $\ell$-adique $b_2(\Xbar)$  de
$\Xbar$ est {\em impair.} En effet, si par malheur ce nombre est pair, on \'eclate un
point ferm\'e de degr\'e $f$ impair (d'apr\`es un argument de type Lang--Weil, un tel
point existe puisque
  la vari\'et\'e $X$ est g\'eom\'etriquement int\`egre),
   ce qui donne pour l'\'eclat\'e $X^*$ la formule $b_2(X^*)=b_2(X)+f$
d'apr\`es \cite{sga7}, expos\'e XVIII, (2.3.1). La raison pour cette hypoth\`ese
suppl\'ementaire se d\'evoilera lors de la preuve de la proposition \ref{cica} ci-dessous.

Un calcul simple de classes de Chern
(voir \cite{schoenams}, 9.2.1) montre que, quitte \`a composer 
le plongement projectif donn\'e de $X$
 avec un plongement de
Veronese de degr\'e {\em pair}, on peut supposer que le deuxi\`eme nombre de Betti
$\ell$-adique de toute section hyperplane lisse de $X$  est  {\em pair}.
Cette information de parit\'e sera
\'egalement importante pour la suite.

Ayant fait une extension de degr\'e premier \`a $\ell$ si n\'ecessaire, on trouve un
\'eclat\'e $V\to X$ muni d'un pinceau de Lefschetz  $V\to D\cong {\bf P}^1$ de sections
hyperplanes (\cite{sga7}, expos\'e XVII, th\'eor\`eme 2.5). Notons $\dot D\subset D$ le
lieu au-dessus duquel le morphisme $V\to D$ est lisse, et choisissons un point
g\'en\'erique g\'eom\'etrique $\varepsilon$ de $D$. D'apr\`es ce qui pr\'ec\`ede,
le deuxi\`eme nombre de Betti de
$V_\varepsilon$ est pair.

Introduisons les notations $\pi$ (resp. $\bar \pi$) pour le groupe fondamental
arithm\'etique (resp. g\'eom\'etrique) de $\dot D$ ayant $\varepsilon$ pour point base.

\begin{prop}\label{imageinfinie}
Quitte \`a faire une extension de $k$ de degr\'e premier \`a $\ell$, on peut choisir le
pinceau $V$ de sorte que l'image de $\bar \pi$ dans  ${\rm Aut}_{\Z_\ell}(H^2(V_\epsilon,
\Z_\ell(1)))$ via la repr\'esentation de monodromie soit infinie.
\end{prop}

\begin{demo} C'est la Proposition 1.1 de \cite{schoen}. On ne donne que l'id\'ee de
l'argument. Soit $\PP$ l'espace projectif param\'etrisant les intersections de $X$ avec
les hypersurfaces de degr\'e fixe suffisamment grand dans un plongement projectif fix\'e.
Soient ${\bf V}\subset \PP\times X$ l'hypersurface universelle, et $\dot \PP\subset \PP$
l'ouvert de lissit\'e de la fibration ${\bf V}\to \PP$. Le choix d'un pinceau de Lefschetz
correspond au choix d'une droite $D\subset\PP$, et on a $V={\bf V}\times_{\PP}D$. Par un
argument de type  Bertini (voir par exemple \cite{fgbook}, Lemma 5.7.2)
apr\`es une extension finie de $k$ de degr\'e premier  \`a $\ell$ on trouve
$D$ assez g\'en\'erale pour laquelle l'homomorphisme $\pi_1(\dot D_\kbar, \epsilon)\to \pi_1(\dot
\PP_\kbar, \epsilon)$ est surjectif. Il suffit donc de montrer
que l'image du deuxi\`eme groupe dans ${\rm Aut}_{\Z_\ell}(H^2(V_\epsilon, \Z_\ell(1)))$
est infinie. Schoen montre par une construction de g\'eom\'etrie alg\'ebrique classique
qu'il existe un autre espace projectif $P$ et un morphisme $P\to\PP$ tels que le morphisme
${\bf V}\times_{\PP}P\to P$ se factorise en ${\bf V}\times_{\PP}P\to W\to P$, o\`u $W$ est
une hypersurface projective lisse, et la dimension relative de $W\to P$ est 2. Comme $W$
est une hypersurface, elle se rel\`eve en caract\'eristique 0, et un th\'eor\`eme de
Deligne (\cite{verdier}, Th\'eor\`eme B) montre que la monodromie de  tout pinceau de
Lefschetz balayant $W$ est infinie (sous l'hypoth\`ese ${\rm car}\,(k)\neq 0$; en
caract\'eristique 2 un petit argument suppl\'ementaire est donn\'e dans \cite{schoen}).
Ceci implique que la monodromie doit \^etre infinie pour la fibration ${\bf
V}\times_{\PP}P\to P$, et finalement pour ${\bf V}\to\PP$.
\end{demo}
\enddem
\smallskip

Expliquons maintenant l'id\'ee de la preuve du th\'eor\`eme \ref{theoschoen}. Tout
d'abord, il suffit de montrer que toute classe dans $H^4(\Xbar, \Z_\ell(2))^G$ est la
classe d'un cycle alg\'ebrique sur $\Xbar$ (ensuite, pour un sous-groupe ouvert $U\subset
G$ on peut appliquer ce r\'esultat apr\`es changement de base de $X$ au sous-corps fix\'e
par $U$). Etant donc donn\'e $w\in H^4(\Xbar, \Z_\ell(2))$ fix\'e par $G$, on montre qu'il
est la poussette d'un \'el\'ement de $H^2(V_{\bar x}, \Z_{\ell}(1))^{\gal(\kbar|k(x))}$,
o\`u $\bar x$ est un point g\'eom\'etrique au-dessus d'un point ferm\'e $x\in\dot D$. La
conjecture de Tate pour les diviseurs sur la surface $V_x$ (qui est valable \`a
coefficients $\Z_\ell$ si elle est valable \`a coefficients $\Q_\ell$ d'apr\`es ce qu'on a
dit dans la remarque \ref{remconjfausse} (3) ci-dessus) montre alors que cet \'el\'ement
est la classe d'un cycle alg\'ebrique.

Notons qu'\`a coefficients $\Q_\ell$ l'\'enonc\'e voulu est une cons\'equence directe du
th\'eor\`eme de Lefschetz difficile ({\em cf.} la preuve du lemme \ref{lemcica} {\em
infra}); toute la finesse de l'argument consiste \`a en tirer un \'enonc\'e \`a
coefficients entiers.

Voici une reformulation. \'Ecrivons $X_\varepsilon$ pour le changement de base
$X\times_{\spec k}\varepsilon$. Par d\'efinition, il est muni de l'action triviale de
$\bar\pi$ (celui-ci agissant sur les fibres de la fibration triviale $X\times D\to D$), et
par cons\'equent
$$H^i(\Xbar, \Z_\ell(2))^{\gal(\kbar|k)}\cong H^i(X_\varepsilon, \Z_\ell(2))^\pi$$ pour
tout $i>0$. Notant $D_{\bar x}$ le groupe de d\'ecomposition d'un point $\bar x$ du
rev\^etement universel profini de $\dot D$ au-dessus de $x$, on obtient
$$H^i(V_{\bar x}, \Z_\ell(2))^{\gal(\kbar|k(x))}\cong
H^i(V_\varepsilon, \Z_\ell(2))^{D_{\bar x}}$$
pour tout $i>0$ par
   le th\'eor\`eme de
changement de base propre et l'isomorphisme $D_{\bar x}\cong{\gal(\kbar|k(x))}$.

Notons $i$ l'inclusion de la surface $V_\varepsilon$ dans la vari\'et\'e $X_\varepsilon$
(qui est de dimension 3). Elle induit un morphisme de restriction
$$ i^*  :\,H^2(X_\varepsilon, \Z_\ell(1))\to H^2(V_\varepsilon, \Z_\ell(1))$$  ainsi qu'une
poussette $$i_*  :\,H^2(V_\varepsilon, \Z_\ell(1))\to H^4(X_\varepsilon, \Z_\ell(2)).$$

D'apr\`es la discussion ci-dessus, il suffit  donc de montrer  :

\begin{prop}\label{cica} Chaque \'el\'ement  de $H^4(X_\varepsilon, \Z_\ell(2))^\pi$ est de la forme
$i_*(\beta)$, avec un $\beta\in H^2(V_\varepsilon, \Z_\ell(1))$ invariant sous l'action
d'un sous-groupe de d\'ecomposition $D_{\bar x}$ dans $\pi$,
pour un point  ${\bar x}$ convenable.
\end{prop}

Interrompons-nous pour  quelques consid\'erations d'alg\`ebre $\Z_\ell$-lin\'eaire.

\section{Le th\'eor\`eme de Schoen, II   : Lemmes d'alg\`ebre lin\'eaire}\label{secalglin}

Etant donn\'es un  $\Z_\ell$-module $B$ et un sous-ensemble $A\subset B$,
 on d\'efinit {\em le satur\'e}
$A_s$ de $A$ dans $B$ comme l'ensemble des $b\in B$ avec $\ell^nb\in A$ pour un $n \geq 0$
convenable. On dit que $A$ est satur\'e dans $B$ si $A_s=A$.

\begin{lem}\label{la1} Pour un module $B$ de type fini sur $\Z_\ell$ il existe un sous-groupe ouvert
$\Gamma\subset Aut_{\Z_l}(B)$
 tel que $B^g$ soit satur\'e dans $B$ pour tout
$g\in \Gamma$.
\end{lem}

\begin{demo} Ecrire $B=F\oplus T$ avec $F$ libre et $T$ de torsion, et prendre
$\Gamma=Aut_{\Z_\ell}(F)\times \{\id_T\}$.
\end{demo}
\enddem

\begin{prop}\label{la2}
Soient $F$ un $\Z_\ell$-module libre de rang fini impair, $S\subset F$ un sous-ensemble
ouvert pour la topologie $\ell$-adique, et $\Phi : F\times F\to\Z_\ell$ une forme
bilin\'eaire sym\'etrique non d\'eg\'en\'er\'ee sur $\Q_\ell$. Il existe alors un
sous-ensemble ouvert ${\mathcal S}\subset O(F, \Phi)$ tel que :

(a) chaque \'el\'ement de
$\mathcal S$ admet un vecteur fixe non nul dans $S$;

(b)  l'ouvert $\mathcal S$ contient des
\'el\'ements arbitrairement proches de $1$ pour la topologie $\ell$-adique.
\end{prop}

Ici $O(F, \Phi)$ d\'esigne le groupe des automorphismes  $\Z_\ell$-lin\'eaires de $F$
pr\'eservant $\Phi$. Pour la preuve nous avons besoin d'un r\'esultat ancillaire.

\begin{lem}\label{lemla1}
Soient $K$ un corps de caract\'eristique diff\'erente de 2, $V$ un \mbox{$K$-espace}
vectoriel de dimension finie impaire, et $\Phi$ une forme quadratique  non
d\'eg\'en\'er\'ee sur $V$. Tout \'el\'ement de $SO(\Phi)$ admet $1$ comme valeur propre.
\end{lem}

\begin{demo} Notons $A$ la matrice de
$\Phi$ dans une base fix\'ee de $V$. Pour un \'el\'ement de $O(\Phi)$ dont la matrice est
$M$ et la matrice transpos\'ee $M^t$, on a
  $M^t .A .M=A$. D'o\`u
  $$M^t .A .(M-I)=A- M^t .A= (I-M^t).A.$$
  En prenant le d\'eterminant on obtient  :
  $$ \det(M). \det(A). \det(M-I)= \det(I-M). \det(A).$$
  Ici $\det(A)\neq 0$ et $\det(M)=1$ (comme $M \in SO(\Phi)$), donc
  $\det(M-I)=\det(I-M)$. Comme $V$ est de dimension impaire, ceci n'est possible que si $\det(M-I)=0$.
\end{demo}
\enddem

\bigskip

\noindent{\em D\'emonstration de la proposition \ref{la2}.} Soit $U\subset
SO(F_{\Q_\ell},\Phi)$ l'ouvert de Zariski form\'e des \'el\'ements ayant des valeurs
propres distinctes. C'est aussi un ouvert de Zariski de $O(F_{\Q_\ell},\Phi)$. Comme $F$
est de rang impair, tout \'el\'ement de $SO(F_{\Q_\ell}, \Phi)$ admet 1 comme valeur
propre par le  lemme \ref{lemla1}. Donc tout $u\in U$ stabilise un sous-espace $L_u$ de
dimension~1 correspondant \`a la valeur propre 1.
 Envoyant $u$ sur $L_u$ on obtient une application continue
$\lambda  :\,U\to {\bf P}(F_{\Q_\ell})$. L'image
de $S \setminus 0$ dans ${\bf P}(F_{\Q_\ell})$
par la projection naturelle $F_{\Q_\ell} \setminus 0 \to {\bf P}(F_{\Q_\ell})$
est
ouverte, tout comme son image inverse $\mathcal S$ dans $U\subset SO(F_{\Q_\ell}, \Phi)$.

Reste \`a voir que l'ensemble $\mathcal S$ est non vide, et qu'il contient des
\'el\'ements arbitrairement proches de~1. Soit $v\in S$ un vecteur  non isotrope. Ecrivant
$F_{\Q_\ell}={\langle v\rangle\perp M}$ avec un $\Q_\ell$-vectoriel $M$, on commence par
montrer qu'il existe un \'el\'ement de $SO(M, \Phi|_M)$ \`a valeurs propres distinctes, et
toutes diff\'erentes de~1.
 Pour cela, on  d\'ecompose l'espace quadratique $M$ en une somme orthogonale d'espaces quadratiques
 $V_i$ de dimension~2. Chaque $SO(V_{i})$ est un tore
$T_{i}=R^1_{k_{i}/k}\G$ de dimension~1, o\`u $k_{i}/k$ est une alg\`ebre \'etale de
degr\'e 2 sur $\Q_\ell$. Si $k_{i} \simeq \Q_\ell \times \Q_\ell$, alors $T_{i} \simeq
{\bf G}_{m,k}$, et ${\bf G}_{m,k}$ agit sur $V_{i}=\Q_\ell \oplus \Q_\ell$ par
$\lambda.(u,v)=(\lambda.u, \lambda^{-1}.v)$. Si $k_{i}$ est une extension quadratique de
$\Q_\ell$, alors $SO(V_{i})(\Q_\ell)$ est le groupe des \'el\'ements de norme 1 dans
$k_{i}$, et l'action de ce groupe sur $V_{i} \simeq k_{i}$ est donn\'ee par la
multiplication dans $k_{i}$.
   Les deux valeurs propres  d'un \'el\'ement $\alpha \in  SO(V_{i})(\Q_\ell) \subset
SO(\Phi) \subset GL(V) $ sont les conjugu\'es de $\alpha$ (qui sont inverses l'un \`a
l'autre).
On trouve donc une famille d'\'el\'ements $\alpha_i \in  SO(V_{i})$  dont la somme d\'efinit un \'el\'ement
de $SO(M,\Phi|_M)$ qui
a des valeurs propres distinctes et diff\'erentes de 1. De plus, on peut
choisir les matrices des $\alpha_i$ de sorte qu'elles aient des
coefficients
dans $\Z_\ell$
et qu'elles soient arbitrairement proches de la matrice 1 pour la topologie $\ell$-adique.
Si elles sont suffisamment
proches de 1, leur somme directe doit pr\'eserver la trace du r\'eseau $F$ sur $M$.
\enddem

\section{Le th\'eor\`eme de Schoen, III   : fin de la d\'emonstration}\label{schoenfin}

Il nous reste \`a d\'emontrer la proposition \ref{cica}.

\begin{lem}\label{lemcica}
Il existe une inclusion
$$
H^4(X_\varepsilon, \Z_\ell(2))^\pi\subset i_*((\ker i_* +H^\pi)_s),
$$
o\`u
$$
H  :=\im(i^*)\subset H^2(V_\varepsilon, \Z_\ell(1)).
$$
\end{lem}

\begin{demo}
Le morphisme compos\'e
$$
H^2(X_\varepsilon, \Z_\ell(1)) \stackrel{i^*}\to H^2(V_\varepsilon,
\Z_\ell(1))\stackrel{i_*}\to H^4(X_\varepsilon, \Z_\ell(2))
$$
est un isomorphisme apr\`es tensorisation par $\Q_\ell$ selon le th\'eor\`eme de Lefschetz
difficile, car c'est le  cup-produit par la classe de la section hyperplane
$V_\varepsilon$. Il  est \'equivariant pour l'action de $\pi$, car $V_\varepsilon$
provient par changement de corps de base d'une $k(\PP^1)$-vari\'et\'e.

 Donc par d\'efinition de $H$
$$i_*(H^\pi)\otimes\Q_\ell\cong
H^4(X_\varepsilon, \Z_\ell(2))^\pi\otimes \Q_\ell,
$$
d'o\`u
\begin{equation}\label{hl}
H^4(X_\varepsilon, \Z_\ell(2))^\pi\subset (i_*(H^\pi))_s. \end{equation}

Remarquons maintenant que le morphisme
$$
i_*  :\,H^2(V_\epsilon, \Z_\ell(1))\to H^4(X_\varepsilon, \Z_\ell(2))
$$
est surjectif. Ceci r\'esulte du th\'eor\`eme de Lefschetz faible  : dans la suite de
localisation
$$
H^4_{V_\varepsilon}(X_\varepsilon, \Z_\ell(2))\to H^4(X_\varepsilon, \Z_\ell(2))\to
H^4(X_\varepsilon\setminus V_\varepsilon, \Z_\ell(2))
$$
le dernier terme est 0, car la vari\'et\'e  $X_\varepsilon\setminus V_\varepsilon$
est affine de dimension 3, et le premier terme est isomorphe
\`a  $H^2(V_\varepsilon, \Z_\ell(1))$ par puret\'e.

En particulier, \'etant donn\'e $w\in H^4(X_\varepsilon, \Z_\ell(2))^\pi$, on trouve
$\beta\in H^2(V_\varepsilon, \Z_\ell(1))$ avec
$$
w=i_*(\beta).
$$
D'autre part, (\ref{hl}) implique
$$
\ell^nw=i_*(\gamma)$$ pour un $\gamma\in H^\pi$ convenable et $n \geq 0$. Mais comme
$\ell^nw=i_*\ell^n\beta$, on obtient $i_*(\gamma-\ell^n\beta)=0$, i.e. $\ell^n\beta\in
\ker i_*+ H^\pi$, d'o\`u le lemme.
\end{demo}
\enddem

\begin{cor} Pour $w\in H^4(X_\varepsilon, \Z_\ell(2))^\pi$ fix\'e, le sous-ensemble
$$
H_w  :=\{v\in \ker i_*  :\, w\in i_*((v+H^\pi)_s)\}
$$
de $\ker i_* \subset H^2(V_\varepsilon, \Z_\ell(1))$ est un ouvert non vide de $\ker i_*$,
stable par multiplication par $\ell$.
\end{cor}

\begin{demo}
Le lemme donne $H_w\neq\emptyset$; plus pr\'ecis\'ement, la preuve du lemme montre que
$v_0  :=\ell^n\beta - \gamma   \in H_w$. Ce choix de $n$ donne $(v_0+\ell^n\ker i_*)\subset
H_w$, car pour $\delta\in\ker i_*$ et $v=v_0+\ell^n\delta$ on a $i_*(\beta+\delta)=w$ et
$\ell^n(\beta+\delta)=v_0+\ell^n\delta+\gamma=v+\gamma\in (v+H^\pi)$. Enfin, la
stabilit\'e de $H_w$ par multiplication par $\ell$ r\'esulte de la d\'efinition.
\end{demo}
\enddem

Consid\'erons maintenant la forme $\Z_\ell$-bilin\'eaire sur $H^2(V_\varepsilon,
\Z_\ell(1))$ induite par le cup-produit (i.e. la forme d'intersection)
sur la cohomologie de la
surface $V_\varepsilon$,
et notons
$H^\perp$ l'orthogonal de $H$. On a alors $\ker i_*\subset H^\perp
 \subset H^2(V_\varepsilon, \Z_\ell(1)),$
 et   l'inclusion ${\ker i_*\subset H^\perp}$ devient \'egalit\'e
apr\`es tensorisation par $\Q_\ell$. En effet, les accouplements de cup-produit satisfont
\`a la compatibilit\'e
$$
\alpha\cup i_*(\beta)=i^*(\alpha)\cup\beta
$$
pour $\alpha\in H^2(X_\varepsilon, \Z_\ell(2))$ et $\beta\in H^2(V_\varepsilon,
\Z_\ell(1))$, et ils sont non d\'eg\'en\'er\'es \`a coefficients $\Q_\ell$.

Soit $F\subset H^\perp$ un $\Z_\ell$-module libre, compl\'ement direct au sous-module de
torsion $T$. Alors $F\cap \ker i_*$ est un sous-module ouvert dans $F$, et d'indice fini
dans $\ker i_*$. Ainsi le corollaire pr\'ec\'edent implique :

\begin{cor}\label{ouvertSw} Pour $w\in H^4(X_\varepsilon, \Z_\ell(2))^\pi$ fix\'e le sous-ensemble
$$
S_w  :=\{v\in \ker i_*\cap F  :\, w\in i_*((v+H^\pi)_s)\}
$$
est un ouvert non vide de $F$.
\enddem
\end{cor}

\begin{rema}\rm
Quand $X$ est une hypersurface dans ${\bf P}^4$, tous les \mbox{$\Z_\ell$-modules}
consid\'er\'es sont sans torsion, et l'on a $\ker i_*=H^\perp=F$, d'o\`u $H_w=S_w$.
\end{rema}

Le lemme suivant distille la strat\'egie de la d\'emonstration de la proposition
\ref{cica}.

\begin{lem}
Fixons $w\in H^4(X_\varepsilon, \Z_\ell(2))^\pi$. Supposons qu'il existe $g\in \pi$
satisfaisant aux trois hypoth\`eses suivantes  :
\begin{enumerate}
\item $H^2(V_\varepsilon, \Z_\ell(1))^g$ est satur\'e dans $H^2(V_\varepsilon, \Z_\ell(1))$;
\item $g$
engendre
topologiquement le
 sous-groupe de d\'ecomposition  $D_{\bar x}$ dans $\pi$;
\item $g$ fixe un \'el\'ement $v\in S_w$.
\end{enumerate}
Alors il existe $\beta\in H^2(V_\varepsilon, \Z_\ell(1))^{D_{\bar x}}$ avec
$w=i_*(\beta)$.
\end{lem}

\begin{demo} Pour un $v$ comme dans (3) on a $(v+H^\pi)\subset H^2(V_\varepsilon, \Z_\ell(1))^g$.
Comme $H^2(V_\varepsilon, \Z_\ell(1))^g$ est satur\'e dans $H^2(V_\varepsilon,
\Z_\ell(1))$, on a de plus $(v+H^\pi)_s\subset H^2(V_\varepsilon,
\Z_\ell(1))^g=H^2(V_\varepsilon, \Z_\ell(1))^{D_{\bar x}}$. Mais par le corollaire
pr\'ec\'edent on a $w=i_*(\beta)$ pour un $\beta\in(v+H^\pi)_s$.
\end{demo}
\enddem

\bigskip

\noindent {\em D\'emonstration de la proposition \ref{cica}.} On cherche un  $g\in \pi$
satisfaisant aux conditions du lemme.

   Par le th\'eor\`eme de Lefschetz
difficile,   la res\-triction de la forme d'inter\-section sur $H^2(V_\varepsilon, \Z_\ell(1))$ \`a
${H\otimes\Q_\ell}$ est non d\'eg\'en\'er\'ee (voir \cite{deligne}, Lemme 4.1.2).
 Sa restriction \`a
$H^\perp\otimes\Q_\ell$ est donc   non d\'eg\'en\'er\'ee. Ecrivant $H^\perp=F\oplus T$
comme ci-dessus, on peut identifier $O(F)$ avec le stabilisateur (point par point) de $T$,
qui est un sous-groupe ouvert d'indice fini de $O(H^\perp)$. Comme l'image de $\bar \pi$
par la repr\'esentation de monodromie $\rho$ est infinie par construction (Proposition
\ref{imageinfinie}),
 un th\'eor\`eme de Deligne (\cite{deligne}, Th\'eor\`eme 4.4.1)
assure que c'est un sous-groupe
ouvert de $O(H^\perp\otimes\Q_\ell)$. A fortiori
$\rho(\pi)\cap O(F)$ est ouvert dans $O(F)$. Par le lemme \ref{la1} il existe un
sous-groupe ouvert $G_0\subset \rho(\pi)\cap O(F)$ tel que $H^2(V_\varepsilon,
\Z_\ell(1))^g$ est satur\'e dans $H^2(V_\varepsilon, \Z_\ell(1))$ pour tout $g\in G_0$.

On applique maintenant la proposition  \ref{la2} \`a $F$. Pour ce faire, on doit d'abord
v\'erifier que le rang de $F$ est impair,  i.e. que la dimension de  $H^\perp\otimes\Q_\ell$ est
impair. Ceci r\'esulte de nos hypoth\`eses initiales que  la dimension de
$H^2(V_\varepsilon, \Q_\ell(1))$ est paire et celle de $H^2(X_\varepsilon, \Q_\ell(1))$
impaire; on  conclut par l'injectivit\'e de $i^*\otimes\Q_\ell$ (voir le d\'ebut de
la preuve du lemme \ref{lemcica}).
La proposition  \ref{la2} (a)
fournit donc un ouvert $\mathcal S$ de $O(F)$
dont tout \'el\'ement a un vecteur fixe dans l'ouvert non vide $S_w$
donn\'e par
le corollaire \ref{ouvertSw}.
 De plus, la proposition \ref{la2} (b) assure que $\mathcal S$ contient des
\'el\'ements arbitrairement proches de 1, donc son intersection avec le sous-groupe ouvert
$G_0$ est un ouvert non vide.

L'image inverse de ${\mathcal S}\cap G_0$ dans $\pi$ est un ouvert dont les \'el\'ements
satisfont aux  propri\'et\'es (1) et (3) du lemme ci-dessus. En outre, par d\'efinition de
la topologie de $\pi$ elle contient une cosette $hV$ d'un sous-groupe normal ouvert
$V\subset \pi$. Appliquant le th\'eor\`eme de densit\'e de Tchebotarev au rev\^etement
galoisien $Z\to \dot D $ correspondant \`a $V$ on obtient un point ferm\'e $z\in Z$ dont le
Frobenius associ\'e dans $\pi/V$ est $\bar h$, la classe de $hV$. Prenons alors un point
$\bar x$ du rev\^etement universel profini de $\dot D$ au-dessus de $z$. Le sous-groupe de
d\'ecomposition $D_{\bar x}$ est engendr\'e par un \'el\'ement $g$ d'image $\bar
h$ dans $\pi/V$. Ceci veut dire que $g$ est un \'el\'ement de $hV$, et en tant que tel
satisfait aux hypoth\`eses (1) et (3) du lemme. Par construction, il satisfait \'egalement
\`a (2).
\enddem

\begin{rema}\rm Si l'on pouvait choisir $V$ de telle  sorte que le
 conoyau du morphisme $\bar \pi/ (V\cap\bar \pi)\to\pi/V$ soit d'ordre premier \`a
$\ell$, alors une variante plus pr\'ecise de l'argument de Tchebotarev ci-dessus donnerait un point
ferm\'e $x$ de degr\'e premier \`a $\ell$. L'existence d'un tel $V$ impliquerait donc la
conjecture de Tate \`a coefficients $\Z_\ell$ pour les 1-cycles  sur $X$ (en supposant la
conjecture connue pour les surfaces).
\end{rema}

\section{Cons\'equences du th\'eor\`eme de Schoen}\label{brauermanin}

Nous donnons ici des applications du th\'eor\`eme \ref{theoschoen} \`a l'existence de
z\'ero-cycles de degr\'e premier \`a la caract\'eristique sur certaines vari\'et\'es
d\'efinies sur le corps des fonctions d'une courbe au-dessus de la cl\^oture alg\'ebrique
d'un corps fini. Il s'agit de deux \'enonc\'es apparent\'es mais non \'equivalents dont
chacun implique le corollaire \ref{corschoen}.

  \begin{theo}\label{edmonton}
Soit $\kbar$ la cl\^oture alg\'ebrique d'un corps fini $k$ de carac\-t\'eristique $p$, et
soit $\overline C$ une $\kbar$-courbe propre lisse connexe, de corps des fonctions $F=
\kbar(\overline C).$ Fixons une cl\^oture s\'eparable $\overline F$ de $F$.

Soit $\Xbar$ une $\kbar$-vari\'et\'e projective, lisse, connexe, admettant un
$\kbar$-morphisme projectif et dominant $f   :  \Xbar \to\overline C$ dont la fibre g\'en\'erique
$\Xbar_F $ est lisse et g\'eom\'etriquement int\`egre.

Supposons   :

(i) Le groupe de Picard $\pic {\Xbar_{\overline F}} $ est sans torsion.

(ii) Pour tout premier $\ell$ diff\'erent de $p$,
la partie $\ell$-primaire
du groupe de Brauer $\br \Xbar \subset \br \Xbar_F$ est
finie.

(iii) La $F$-vari\'et\'e $\Xbar_F$ a des points dans tous les compl\'et\'es de $F$ aux
points \mbox{de $C$.}

(iv) Pour tout premier $\ell$ diff\'erent de $p$, la conjecture de Tate $\ell$-adique vaut
pour les diviseurs sur les surfaces projectives et lisses sur un corps fini de
caract\'eristique $p$.

Alors $\Xbar_F$ poss\`ede un z\'ero-cycle de degr\'e une puissance de $p$.
\end{theo}

 \begin{demo}
 Soit $d+1$ la dimension de $\Xbar$, et donc $d$ la dimension de la $F$-vari\'et\'e $\Xbar_F$.
Fixons une cl\^oture s\'eparable ${\overline F}$ de $F$.  Consid\'erons la suite d'applications
$$
CH^{d}(\Xbar) \otimes \Z_\ell \to  H^{2d}(\Xbar, \Z_\ell(d)) \to
H^{2d}(\Xbar_F,\Z_\ell(d)) \to H^{2d}(\overline{X}_{\overline F}, \Z_\ell(d))\cong\Z_\ell
$$
La fl\`eche $H^{2d}(\Xbar, \Z_\ell(d)) \to
H^{2d}(\Xbar_F,\Z_\ell(d))$ est obtenue par passage \`a la limite
sur les applications de restriction $$H^{2d}(\Xbar, \Z_\ell(d)) \to
H^{2d}(\Xbar\times_{\overline C}U, \Z_\ell(d)) $$
pour $U$ parcourant les ouverts non vides de la courbe ${\overline C}$.
La compatibilit\'e de l'application cycle \`a la restriction \`a un ouvert
(\cite{milneEC}, \S VI, Prop.~9.2)
montre que l'application
compos\'ee ci-dessus se factorise \`a travers le groupe ${CH^{d}(\Xbar_{F})\otimes\Z_\ell}$.
 Il suffit
donc d'\'etablir la surjectivit\'e de l'application compos\'ee en question, pour tout
$\ell$ premier \`a $p$. On le fait en plusieurs \'etapes.\smallskip

\noindent {\em Surjectivit\'e de $CH^{d}(\Xbar) \otimes \Z_\ell \to  H^{2d}(\Xbar,
\Z_\ell(d))$.} Il existe un corps fini $k\subset\kbar$ et une $k$-vari\'et\'e $X$ telle
que $\Xbar=X\times_k\kbar$. La surjectivit\'e requise r\'esulte du  th\'eor\`eme
\ref{theoschoen}, pourvu qu'on d\'emontre que tout \'el\'ement de $H^{2d}(\Xbar,
\Z_\ell(d))$ est fix\'e par un sous-groupe ouvert de $\gal(\kbar|k)$.

Or pour tout $n>0$ la dualit\'e de Poincar\'e $$ H^{2d}(\Xbar,\mu_{\ell^n}^{\otimes d}))
\times H^2(\Xbar,\mu_{\ell^n}) \to \Z/\ell^n\Z$$ est un accouplement parfait \'equivariant
pour l'action de Galois. D'autre part, on a la suite exacte de Kummer
$$0 \to \pic \Xbar/\ell^n\pic\Xbar  \to H^2(\Xbar,\mu_{\ell^n}) \to {}_{\ell^n}\br \Xbar \to 0.$$
Le groupe $\pic \Xbar$ est extension du
groupe de N\'eron-Severi $NS(\Xbar)$ par le groupe $\ell$-divisible $\pic^0\,\Xbar$. Le
groupe $NS(\Xbar)$ est de type fini, donc $NS(\Xbar)/\ell^n\cong \pic\Xbar/\ell^n$ est fix\'e
par un sous-groupe ouvert de $\gal(\kbar|k)$ ind\'ependant de $n$. Par l'hypoth\`ese
$(ii)$ le groupe ${}_{\ell^n}\br \Xbar$ a \'egalement cette propri\'et\'e. Il en est donc
de m\^eme pour le groupe fini $H^{2d}(\Xbar,\mu_{\ell^n}^{\otimes d})$, et a fortiori pour
$H^{2d}(\Xbar,\Z_\ell(d))$.\smallskip

\noindent{\em Surjectivit\'e de $H^{2d}(\Xbar, \Z_\ell(d)) \to
H^{2d}(\Xbar_F,\Z_\ell(d))$.} On va d\'emontrer la surjectivit\'e de $H^{2d}(\overline
X,\mu_{\ell^n}^{\otimes d}) \to H^{2d}(\overline{X}_{F}, \mu_{\ell^n}^{\otimes d})$ pour
tout $n>0$. Ceci donnera un morphisme surjectif de syst\`emes projectifs de groupes
ab\'eliens finis, d'o\`u une surjection apr\`es passage \`a la limite projective suivant
$n$.

 On a la suite exacte de localisation
$$H^{2d}(\Xbar,\mu_{l^n}^{\otimes d}) \to H^{2d}(\Xbar_F, \mu_{l^n}^{\otimes d}) \to
 \bigoplus_{P\in \overline C_0} H_{\Xbar_{P}}^{2d+1}(\Xbar, \mu_{l^n}^{\otimes d}),$$
 o\`u $P$ parcourt les points ferm\'es de $\overline C$.
 Montrons que chaque fl\`eche
$H^{2d}(\Xbar_F, \mu_{l^n}^{\otimes d}) \to
  H_{\Xbar_{P}}^{2d+1}(\Xbar, \mu_{l^n}^{\otimes d})$ est nulle.  Pour ce faire, par excision
 on peut se restreindre au-dessus du hens\'elis\'e
 $R=O_{\overline C,P}^h$  de $\overline C$ en $P$, dont on note $L$ le corps des fractions.  On dispose de la suite exacte de localisation
 $$ H^{2d}(\overline X_{R}, \mu_{l^n}^{\otimes d}) \to H^{2d}(\overline X_{L}, \mu_{l^n}^{\otimes d})
 \to H_{X_{P}}^{2d+1}(\overline X_R, \mu_{l^n}^{\otimes d}).$$
Il suffit donc d'\'etablir la surjectivit\'e du morphisme $H^{2d}(\Xbar_{R},
\mu_{\ell^n}^{\otimes d}) \to H^{2d}(\Xbar_{L}, \mu_{\ell^n}^{\otimes d})$.

Le corps $L$ est de dimension cohomologique 1   (c'est un corps $C_{1}$).
 La suite spectrale de Hochschild--Serre donne donc naissance \`a la suite exacte courte
{\small $$
 0 \to H^1(L,H^{2d-1}(\overline{X}_{\overline L}, \mu_{\ell^n}^{\otimes d}))  \to H^{2d}(\overline
X_L,\mu_{\ell^n}^{\otimes d})\to H^{2d}(\overline{X}_{\overline L}, \mu_{\ell^n}^{\otimes
d})^{\gal(\overline L|L)} \to 0.
$$}

\noindent D'autre part, le groupe $H^{2d-1}(\overline{X}_{\overline L},
\mu_{\ell^n}^{\otimes
 d})$ est dual de $H^{1}(\overline{X}_{\overline L}, \mu_{\ell^n})\cong
 {}_{\ell^n}\pic\Xbar_{\overline L}$ par
 dualit\'e de Poincar\'e. La torsion dans le groupe de Picard ne change pas
 par extension de corps alg\'ebriquement clos, donc ce dernier groupe est nul par
 l'hypoth\`ese $(i)$. Ceci donne des isomorphismes de modules galoisiens
 $$ H^{2d}(\Xbar_{L}, \mu_{\ell^n}^{\otimes d})) \cong H^{2d}( \Xbar_{\overline{L}},
  \mu_{\ell^n}^{\otimes d})\cong\Z/\ell^n\Z.$$
L'hypoth\`ese que $\Xbar_F$ poss\`ede des points dans tous les compl\'et\'es de $F$ est
\'equivalente \`a la m\^eme hypoth\`ese avec les hens\'elis\'es, d'o\`u en particulier une
section du morphisme  $\Xbar_{R} \to \spec R$ par propret\'e de $\Xbar$. Elle donne
naissance \`a un 1-cycle sur $\Xbar_{R}$ dont la classe de cohomologie dans
$H^{2d}(\Xbar_{R},\mu_{\ell^n}^{\otimes d})$ s'envoie sur le g\'en\'erateur de  $
H^{2d}(\Xbar_{L}, \mu_{\ell^n}^{\otimes d})$.

\smallskip

\noindent{\em Surjectivit\'e de $H^{2d}(\overline X_F,\Z_\ell(d)) \to
H^{2d}(\overline{X}_{\overline F}, \Z_\ell(d))$.}  Ici encore, il suffit d'\'etablir la
surjectivit\'e \`a niveau fini, car il r\'esulte de l'\'etape pr\'ec\'edente que les
groupes $H^{2d}(\overline{X}_{F}, \mu_{\ell^n}^{\otimes d})$ sont finis. Le corps $F$ est
de dimension cohomologique 1, donc le morphisme
$$
H^{2d}(\overline X_F,\mu_{\ell^n}^{\otimes d}) \to H^{2d}(\overline{X}_{\overline F},
\mu_{\ell^n}^{\otimes d})^{\gal(\overline F|F)}$$
dans  la suite spectrale de Hochschild--Serre
est une surjection.
On a
$H^{2d}(\overline{X}_{\overline F}, \mu_{\ell^n}^{\otimes d})\cong \Z/\ell^n\Z$ avec
action triviale de Galois, d'o\`u la surjectivit\'e requise.
 \end{demo}\
 \enddem

\medskip

Une cons\'equence facile du th\'eor\`eme est le corollaire \ref{corschoen} de
l'intro\-duction.\medskip

\noindent{\em D\'emonstration du corollaire \ref{corschoen}.\/} Comme le degr\'e de la
fibre g\'en\'erique du morphisme $\Xbar\to \overline C$ est suppos\'e premier \`a $p$, il
suffit de montrer que les hypoth\`eses $(i)$ et $(ii)$ du th\'eor\`eme \ref{edmonton} sont
satisfaites. En d'autre termes, on doit v\'erifier que pour $\overline X_{\overline F}$
une intersection compl\`ete lisse de dimension $\geq 3$ dans $\PP^n_{\overline F}$ le
groupe de Picard est sans torsion $\ell$-primaire et la partie $\ell$-primaire du groupe
de Brauer est finie.

Le premier \'enonc\'e r\'esulte du th\'eor\`eme de Noether--Lefschetz (\cite{sga2},
expos\'e XII, corollaire 3.7): la fl\`eche de restriction $\Z= \pic \PP^n_{\overline F}
\to \pic \overline X_{\overline F}$ est un isomorphisme. D'autre part, la suite exacte de
Kummer en cohomologie \'etale donne une suite exacte
$$
0\to \pic \Xbar_{\overline F}/\ell\,\pic \Xbar_{\overline F}\to H^2(\Xbar_{\overline F},
\Z/\ell\Z)\to{}_\ell\br \Xbar_{\overline F}\to 0.
$$
On vient de voir que le premier terme est isomorphe \`a $\Z/\ell\Z$. Mais ceci vaut
\'egalement pour le deuxi\`eme, car il est isomorphe \`a $H^2(\PP^n_{\overline
F},\Z/\ell\Z)$ par le th\'eor\`eme de Lefschetz faible en cohomologie \'etale. On constate
donc avec satisfaction que le dernier terme disparait, ce qui montre que la partie
$\ell$-primaire de $\br \Xbar_{\overline F}$ est en fait triviale.
\enddem

Le corollaire \ref{corschoen} peut \'egalement se d\'eduire du th\'eor\`eme
\ref{vraibrauermanin} ci-apr\`es. Il s'agit d'une variante du th\'eor\`eme  \ref{edmonton}, avec la
diff\'erence sensible qu'ici on fait une hypoth\`ese au-dessus du corps de fonctions d'une
courbe d\'efinie sur un corps fini, et non pas sur $\overline\F_p$.

Soient donc $C$ une courbe propre, lisse, g\'eom\'etriquement connexe  d\'efinie sur un
corps fini $k$, et $Y$ une vari\'et\'e lisse sur le corps des fonctions $k(C)$ de $C$. Pour un
point ferm\'e $P$ de $C$ notons $K_P$ le compl\'et\'e de $k(C)$ pour la valuation
discr\`ete associ\'ee \`a $P$. Une famille $\{z_P\}$ de z\'ero-cycles de degr\'e 1 sur
$Y\times_{k(C)}K_P$ pour tout $P$ d\'efinit un homomorphisme $\br Y\to \Q/\Z$ donn\'e par
$A\mapsto \Sigma_P{\rm inv}_P (A[z_P])$. Ici $\br Y$ est le groupe de Brauer de $Y$, le
morphisme ${\rm inv}_P$  est l'invariant de Hasse du corps local $K_P$, et $A[z_P]$ est
l'\'evaluation de l'\'el\'ement $A\in\br Y$ en $z_P$ d\'efinie via la fonctorialit\'e
contravariante du groupe de Brauer.

On dit qu'il n'y a pas d'obstruction de Brauer--Manin \`a l'existence d'un z\'ero-cycle de
degr\'e 1 sur $Y$ s'il existe une famille $\{z_P\}$ pour \mbox{laquelle} l'homomorphisme
ci-dessus est nul. Notons que cette condition est automatique si la fl\`eche naturelle $\br k(C) \to \br Y$
est surjective.

\begin{theo}\label{vraibrauermanin}
     Soient $k$ un corps fini de caract\'eristique $p$, et $X\to C$ un morphisme projectif et dominant
     de $k$-vari\'et\'es projectives
     lisses g\'eom\'etriquement connexes, o\`u $C$ est une courbe et la fibre g\'en\'erique $X_{k(C)}$ est lisse
     et g\'eom\'etrique\-ment int\`egre.

Supposons   :

(i) Il n'y a pas d'obstruction de Brauer--Manin \`a l'existence d'un z\'ero-cycle de
degr\'e 1 sur la $k(C)$-vari\'et\'e~$X_{k(C)}$.

(ii) La conjecture de Tate sur les diviseurs vaut sur les surfaces projectives et lisses
sur un corps fini.

Alors    la  $\overline k(C)$-vari\'et\'e
$X\times_{k(C)}\overline k(C)$ poss\`ede un z\'ero-cycle de degr\'e une puissance de $p$.
\end{theo}

\begin{demo}
 Soit $d+1$ la dimension de $X$, et donc $d$ la dimension de la $k(C)$-vari\'et\'e $X_{k(C)}$.
Soit $\ell \neq p$ un nombre premier, et soit $\{z_P\}$ une famille de z\'ero-cycles de
degr\'e 1 sur les $X\times_{k(C)}K_{P}$ pour laquelle l'application $\br X_{k(C)}\to\Q/\Z$
d\'efinie ci-dessus est nulle. La proposition 3.1 de \cite{CT} fournit alors un
\'el\'ement $z$ de $H^{2d}(X,\Z_\ell(d))$ dont la restriction dans le groupe
$H^{2d}(X\times_{k(C)}{\overline{k(C)}},\Z_\ell(d)) \simeq \Z_\ell$ est \'egale \`a $1 \in
\Z_\ell$. D'apr\`es le th\'eor\`eme \ref{theoschoen}, l'image de $z$ dans
${H^{2d}(X\times_k\kbar, \Z_{\ell}(d))}$ provient d'un \'el\'ement $Z$ de
$CH_{1}(X\times_k\kbar) {\otimes} \Z_\ell$. Prenant la trace de $Z$ dans le groupe
${CH_0(X\times_{k(C)}\overline k(C))\otimes\Z_\ell}$ on voit que le morphisme
$CH_0(X\times_{k(C)}\overline k(C))\to\Z/\ell\Z$ induit par le degr\'e est surjectif.
\end{demo}
\enddem

\begin{rema}\rm
La d\'emonstration de la proposition 3.1 de \cite{CT} invoqu\'ee ci-dessus utilise des
arguments voisins de ceux rencontr\'es dans la preuve du th\'eor\`eme \ref{edmonton}. Une
diff\'erence notable est que, dans la notation de ladite preuve, l'existence de classes de
cohomologie de degr\'e 1 sur les sch\'emas $X_R$ implique directement l'existence d'une
classe globale de degr\'e 1 sur $X$ gr\^ace \`a l'hypoth\`ese \og arithm\'etique \fg de type
Brauer--Manin, sans imposer une hypoth\`ese g\'eom\'etrique comme l'hypoth\`ese $(i)$ du
th\'eor\`eme \ref{edmonton}.

  \end{rema}

\bigskip

{\noindent\small{\em Remerciements.} Nous remercions chaleureusement Burt Totaro pour nous
avoir communiqu\'e la d\'emonstration de \ref{AH}, et Bruno Kahn pour plusieurs
  discussions \'edifiantes. Une grande partie de ce travail a \'et\'e fait
  lors d'une visite du premier auteur \`a Budapest dans le cadre du programme
  intra-europ\'een BUDALGGEO de l'Institut R\'enyi. Le second auteur remercie
  \'egalement OTKA (projet no. K 61116) pour son soutien.}

\end{document}